\documentclass[12pt,leqno]{article}

\usepackage{amsmath}
\usepackage{amssymb}

\usepackage[french]{babel}

\makeatletter
\@addtoreset{equation}{section}
\makeatother

\newif\ifpdf
\ifx\pdfoutput\undefined \pdffalse % we are not running PDFLaTeX
\else \pdfoutput=1 % we are running PDFLaTeX
\pdftrue \fi

\ifpdf
\usepackage[pdftex]{graphicx}
\usepackage[pdftex]{hyperref}
\else
\usepackage{graphicx}
\fi

\newtheorem{thm}{Theorem}[section]
\newtheorem{defi}{Definition}

\newtheorem{lem}[thm]{Lemma}

\newtheorem{rema}[thm]{Remark}

\def\N{{\mathbb N}}
\def\R{{\mathbb R}}

\def\VE{\varepsilon}
\def\Im{\mathop{\rm Im}\nolimits}

\def\supp{\mathop {\rm supp}\nolimits}

\def\sgn{\mathop{\rm sgn}\nolimits}

\def\v{\vert}

\def\cqfd{\hbox{\vrule\vbox{\hrule\phantom{o}\hrule}\vrule}}
\def\cqfd{\hfill$\blacksquare$}

\def\ds{\displaystyle}

\title{Brief Article}
\author{The Author}

\title{Remark on the Kato smoothing effect for Schr\"odinger equation with superquadratic potentials}
\author{Luc Robbiano and Claude Zuily}
\date{}

\begin{document}

%\tableofcontents (s'il y a une table des matires)

\maketitle

\begin{abstract}
The aim of this note is to extend recent results of Yajima-Zhang \cite{Y-Z1, Y-Z2} on the $\frac{1}{2}$- smoothing effect for Schr\"odinger equation with potential growing at infinity faster than quadratically.
\end{abstract}

\section{Introduction}
The aim of this note is to extend a recent result by Yajima-Zhang \cite{Y-Z1, Y-Z2}. In this paper these authors considered the Hamiltonian $H=-\Delta+V(x)$ where $V$ is a real and $C^\infty$ potential on $\R^n$ satisfying for some $m>2$ and $\langle x \rangle = (1+\vert x\vert^2)^\frac{1}{2}$,
   \begin{equation}\label{eq1.1} 
 \vert\partial^\alpha V(x)\vert \leq C_\alpha \langle x \rangle^{m-\vert\alpha\vert },\  x\in \R^n,\  \alpha \in \N^n,
 \end{equation}
 \begin{equation}\label{eq1.2} 
 \ {\rm for \ large} \  \vert x\vert, \quad V(x)\geq C_1 \vert x \vert^m,\quad  C_1 > 0, 
 \end{equation}
 and they proved the following. For any $T>0$ and $\chi \in C^\infty_0(\R^n) $ one can find $C>0$ such that for all $u_0$ in $L^2 (\R^n)$,
  \begin{equation}\label{eq1.3} 
      \int^T_0\Vert\chi(I-\Delta)^\frac{1}{2m} e^{-itH}u_0\Vert ^2 _{L^2 (\R^n)}dt \leq C\Vert u_0\Vert _{L^2(\R^n)}^2
   \end{equation}
   where $\Delta$ is the flat Laplacian.
   In this note, using the ideas contained in Do\"i \cite{D3} we shall show that one can handle variable coefficients Laplacian with time dependent potentials, one can remove the condition (\ref{eq1.2}), one can replace the cut-off function $\chi$ in (\ref{eq1.3}) by $\langle x\rangle^{-\frac{1+\nu}{2}} $ with any $\nu>0$ and finally that the weight $\langle x\rangle^{-\frac{1}{2}}$ is enough for the tangential derivatives.\\
   When $V=0$ the estimate (\ref{eq1.3}) goes back to Constantin-Saut \cite{C-S}, Sj\"olin \cite{S}, Vega \cite{V}, Yajima \cite{Y} who extended to the Schr\"odinger equation a phenomenon discovered by T. Kato \cite{K} on the KdV equation. Later on their results where extended to the variable coefficients operators by Do\"i in a series of papers \cite {D1,D2,D3,D4} which contained the case $m=2$ of Theorem 1.1 below.\\
   Let us describe more precisely our result. It will be convenient to introduce the H\" ormander's metric
        \begin{equation}\label{eq1.4} 
    g= \frac{dx^2}{\langle x \rangle^2} + \frac{d\xi^2}{\langle \xi\rangle^2}
      \end{equation}
      to which we associate the usual class of symbols $S(M,g)$ if $M$  is a weight. Recall that $q \in S(M,g)$ iff $q \in C^\infty (\R^{2n})$ and
      \begin{displaymath}
      \forall \alpha,\beta \in \N^n \ \exists C_ {{\alpha}{\beta}}>0,\  \vert \partial^\beta_x\partial^\alpha_{\xi} q(x,\xi)\vert \leqq  C_ {{\alpha}{\beta}} M(x,\xi) \langle x\rangle^{-\vert \beta\vert} \langle{\xi}\rangle^{-\vert
      \alpha\vert},\ \forall (x,\xi)  \in  T^*(\R^n)
   \end{displaymath}
   If $T>0$ we shall set
        \begin{equation}\label{eq1.5}
        S_T(M,g) = L^\infty([0,T], S(M,g)).
       \end{equation}
       We shall consider here an operator $P$ of the form
        \begin{equation}\label{eq1.6} 
        P= \sum_{j,k=1}^n (D_j-a_j(t,x))g^{jk}(x)(D_k-a_k(t,x)) + V(t,x)
               \end{equation}
               and we shall denote by $p$ the principal symbol of $P$, namely
                  \begin{equation}\label{eq1.7} 
                  p(x,\xi) = \sum _{j,k=1}^n g^{jk}(x)\xi_j\xi_k.
                     \end{equation}
                      We shall make the following structure and geometrical assumptions.
                      
                   Structure assumptions.
              We shall assume the following, 
                \begin{equation}\label{eq1.8}\left\{
                \begin{array}{l}
                
             (i) \   {\rm the \ coefficients} \ a_j, g^{jk},\  V {\rm are \  real  \ valuedÊ\ for}\  j,k = 1,...,n ,\\
             (ii)\  p\in S(\langle\xi \rangle^2,g)\  {\rm and} \  \nabla g^{jk}(x)=o(\vert x \vert^{-1}), \ \v x \v \to + \infty \  1\leqq j,k   \leqq, n \\
             (iii)\ a_j \in S_T(\langle x\rangle^{\frac{m}{2}},g),\ 1 \leqq jÊ\leqq n,\ V \in S_T(\langle x\rangle^m,g)\ m \geqq 2
             \end{array}\right.
             \end{equation}
               \begin{equation}\label{eq1.9} 
               \exists \delta >0,\  p(x, \xi) \geqq \delta \vert \xi \vert^2,\quad   \forall (x, \xi) \in T^*(\R^n).
                  \end{equation}                                    
            \begin{equation}\label{eq1.10} 
            {\rm For\  any \  fixed\  t\ ÊÊ in\ } Ê [0,T] \ Ê{\rm the  \  operator}\    ÊP \  {\rm is  \  essentially  \ self  \ adjoint  \ on}\  L^2(\R^n)
                     \end{equation}
                     Geometrical assumptions.  Let $\phi_t $  be the bicharacteristic flow of $p$. It is easy to see that under the conditions (\ref{eq1.8}),\ (\ref{eq1.9}) it is defined for all $t\inÊ\R$. Let us set $S^*(\R^n) =\lbrace (x, \xi) \in T^*(\R^n) : p(x, \xi) = 1 \rbrace$. Then we shall assume that, 
        \begin{equation}\label{eq1.11}
        \forall K {\rm compact }\subseteqq  S^*(\R^n) \   \exists t_K >0\  {\rm such\  that} \  \Phi_t(K)\cap K =\emptyset \  ,\quad    \forall t  \geq t_K.
           \end{equation}    
           This is the so-called "non trapping condition" which is equivalent to the fact that if $\Phi_t(x;\xi)=(x(t),(\xi(t))$ then  $\ds \lim_{t \to +\infty} \vert x(t)\vert = +\infty$.
                            
           We shall consider $u \in C^1([0,T], {\cal S}(\R^n) )$ and we set
                  \begin{equation}\label{eq1.12} 
                  f(t)= (D_t+P)u(t)
                     \end{equation}  
         For $s \in \R$ \ let $e_s(x,\xi)=(1+\vert \xi \vert ^2 + \vert x\v^m)^{\frac {s}{2}} $ and $E_s$ be the Weyl quantized pseudo-differential operator with symbol $e_s$.\\
                           Our first result is the following.
                 \begin{thm}
                 Let $T>0$. Let $P$ be defined  by (\ref{eq1.6}) which satisfies (\ref{eq1.8}), (\ref{eq1.9}), (\ref{eq1.10}),(\ref{eq1.11}). Then for any $\nu >0$ one can find $C=C(\nu, T) >0$ such that for any $u \in C^1([0,T], {\cal S}(\R^n))$ and all $t$ in $[0,T]$  we have,
                      \begin{displaymath}
                \Vert u(t)\Vert ^2_{L ^2} +   \int^T_0\Vert \langle x \rangle^{-\frac{1+\nu}{2}} E_{\frac{1}{m}}u(t)\Vert^2_{L^2}dt \leqq \ C\  (
                     \Vert u(0)\Vert ^2_{L^2}+ \int^T_0\Vert \langle x \rangle^{\frac{1+\nu}{2}} E_{-\frac {1}{m}}f(t)\Vert^2_{L^2}dt ).
                       \end{displaymath}
                       Here $L^2=L^2(\R^n)$ and $f(t)$ is defined by (\ref{eq1.12}).
                          \end{thm}
                       Now even when $P$ is the flat Laplacian it is known that the estimate in the above Theorem
                       does not hold with $\nu=0$. However we have the following result. Let us set
    \begin{equation}\label{eq1.13} 
    \ell_{jk}=\frac {x_j\xi_k - x_k \xi_j}{\langle x \rangle \langle \xi \rangle},\quad 1\leqq j,k \leqq n \ ,
      \end{equation}  
      and let us denote by $\ell ^w_{jk}$ its Weyl quantization.
            \begin{thm}
            Let $T>0$. Let $P$ be defined by (\ref{eq1.6})  with real coefficients satisfying (\ref{eq1.9}),
            (\ref{eq1.10}),(\ref{eq1.11}) and
            \begin{equation}\label{eq1.14}\left\{
            \begin{array}{l}
                                   (i)\  g_{jk}=\delta _{jk}+b_{jk}, \  b_{jk}  \in S(\langle x\rangle^{-\sigma_0},g),\    for\    some\   Ê\sigma_0 >0,\\
             
            (ii)\   a_j \in S_T(\langle x \rangle^\frac{m}{2},g),\   V \in S_T(\langle x \rangle^m,g).
            \end{array}\right.
            \end{equation} 
            Then for any $\nu>0$ one can find $C=C(\nu,T)$ such that for any $u \in C^1([0,T], {\cal S}(\R^n))$  and $f(t)=(D_t+P)u(t)$ we have
              \begin{displaymath}
                 \sum_{j,k=1}^n   \int^T_0\Vert \langle x \rangle^{-\frac{1}{2}} E_{\frac{1}{m}}\ell^w_{jk}u(t)\Vert^2_{L^2}dt \leqq \  C\ (\Vert u(0)\Vert ^2_{L^2}+ \int^T_0\Vert \langle x \rangle^{\frac{1+\nu}{2}} E_{-\frac {1}{m}}f(t)\Vert^2_{L^2} dt).
                       \end{displaymath}
                \end{thm}
                Here are some  remarks and examples.
                \begin {rema}
  1)We know that one can find   $\psi \in C^ \infty _0{(\vert x \vert
   <1)}$ and $\phi \in  C^ \infty _0{(\frac{1}2} \leqq  \vert x \vert \leqq 2)$ positive such that
     $\ds \psi(x)+Ê\sum_{j=0}^{+ \infty}\phi (2^{-j}x)=1$,\  for all x in $\R^n$. Let $\ds V=\vert x \vert^m \sum _{j{\rm\ even}} \phi(2^{-j}x)- \vert x \vert^2 \sum_{j {\rm\ odd}} \phi(2^{-j}x)$. Then $V \in S(\langle xÊ\rangle ^m, g)$ and since $V \geqq -\vert x\vert^2$ the operator $P=-Ê\Delta + V$\ is essentially self adjoint on  $C^ \infty _0(\R^n)$. It follows that (\ref{eq1.9}), (\ref{eq1.10}), (\ref{eq1.11})  and  (\ref{eq1.14}) are satisfied,   therefore Theorem 1.1 and 1.2 apply. However the lower bound (\ref{eq1.2}) assumed in \cite {Y-Z2} is not satisfied.\\
     2) Assume that $\ds p(x,\xi)=\vert \xiÊ\vert^2 + \varepsilon \sum_{j,k=1} ^n b_{jk}(x)\xi_j \xi_k$  with $b_{jk} \in S(\langle x\rangle^{-\sigma_0},g)$  for some $ \sigma_0>0$. Then if $\varepsilon$ is small enough the non trapping condition (\ref {eq1.11}) is satisfied.
     \end{rema}
     \section{Proofs of the results}
     Let us consider the symbol\   $\ds a_0(x,\xi)=\frac {x \cdot\xi}{\langle \xi\rangle}$. A straightforward computation shows that under condition (\ref {eq1.8})  $(ii)$  one can find\  $C_0, C_1, R$ positive such that
         \begin{equation}\label{eq2.1} 
         H_{p}a_0(x,\xi) \geqq C_0 \vertÊ\xiÊ\vert - C_1,\  if\  (x, \xi) \in  
         T^* (\R^n)\  and\   \vert x \vert \geqq R.
     \end {equation}\\
     where $H_p$ denotes the Hamiltonian field of the symbol $p$.\\
           \rm{Then we have the following result due to Do\" i [D3].}
     \begin{lem}
     Assume moreover that (\ref{eq1.11}) is satisfied then there exist $a \in S(\langle x \rangle,g)$ and positive constants $C_2,C_3$ such that\\
         (i) $H_pa(x,\xi) \geqq C_2 \vert \xi\vert - C_3,\quad  \forall (x,\xi) \in  T^*(\R^n)$,\\
         (ii) $a(x,\xi) = a_0(x,\xi)$,\quad if $\vert x\vert $\  is large enough.
           \end{lem}
                  The  symbol $a$ is called a global escape function for $p$. Here is the form of this symbol. Let $\chi \in C^\infty _{0}(\R^n)$ be such that $\chi(x)=1$ if\  $\vert x\vert \leqq 1$, $\chi(x)=0$\ if\ $\vert x\vert \geqq 2$ and $0 \leqq \chi \leqq 1$. With $R$ large enough and $M \geqq 2R$ we have,
    \begin{displaymath}
 a(x,\xi) = a_0(x,\xi) +M^ \frac{1}{2}  \chi \big(\frac{x}{M}\big) a_1\big(x, \frac {\xi}{\sqrt {p(x,Ê\xi)}}\big)(1 -\theta(\sqrt{p(x,\xi)})
   \end{displaymath}
       where 
       \begin{displaymath}
       a_1(x,\xi)= - \int_ {0}^{+\infty} \chi\big(\frac {1}{R} \pi(\Phi_t(x,\xi)\big)\ dt
        \end{displaymath}
        and $\pi(\Phi_t(x,\xi))= x(t;x,\xi)\ , \theta(t)= 1$ if $0\leqq t\leqq 1,\  \theta (t)= 0 $\ if $t\geqq 2,\  0Ê\leqq \theta \leqq 1$. Details can be found in [D3].\\
        Proof of Theorem 1.1\\
        Let $\psi \in C^ {\infty}(\R^n)$ be such that $\supp \psi \in [\varepsilon, {+\infty}[ ,\ \psi(t) =1 $ in $[2\varepsilon, {+\infty}[$\  (where $\varepsilon >0 $ is a small constant chosen later on) and $\psi ' (t) \geqq 0$ for $t \in \R$. Following Do\"i [D3] we set,
                   \begin{equation}\label{eq2.2}\left\{
            \begin{array}{l}
        \psi_0(t) = 1 - \psi(t) - \psi(-t) = 1- \psi(\vert t\vert)\\
        \psi_1(t) = \psi(-t) - \psi(t) = - \sgn t\   \psi(\vert t\vert)
          \end{array}\right.
          \end{equation}
      Then $\psi_j \in C^{\infty} (\R), {\rm for}\   j=0,1$ and we have
    \begin{equation}\label{eq2.3}
                   \psi _{0}'(t) = - \sgn t\  \psi '(\vert t\vert)\quad  {\rm and}\quad
      \psi _{1}'(t) = - \psi '(\vert t\vert).
                 \end{equation}
         Let $\chi \in C^{\infty}(\R)$ be such that $\chi(t) = 1$\  if $ t\leqq \frac {1}{2} ,\  \chi(t) = 0$ if $tÊ\geqq 1$
and $ \chi (t) \in [0,1]$. With $a$  given by Lemma 2.1 we set
         \begin{equation}\label{eq2.4}\left\{
            \begin{array}{l}
           \ds  \theta(x, \xi) = \frac {a(x,\xi)}{\langle x\rangle},\quad    (x, \xi) \in T^*(\R^n),\\
           \ds  r(x, \xi) = \frac {\langle x\rangle^ \frac {m}{2}}{\sqrt {p(x, \xi)}},\quad     (x, \xi) \in T^*(\R^n) \setminus  0.
           \end{array}\right.
          \end{equation}   
          Finally we set
        \begin{equation}\label{eq2.5} 
           -\lambda =\Big( \frac {a} {\langle x\rangle} \psi_0(\theta) - \big(M_0 - \langle     a\rangle^{-\nu}\big) \psi_1(\theta)\Big)  p^{\frac{1}{m}-\frac {1}{2}}  \chi(r),
            \end{equation}
           where $\nu>0$ is an arbitrary small constant and $M_0$ a large constant to be chosen.\\
           The main step of the proof is the following Lemma.
           \begin{lem}
           (i) One can find $M_0>0$ such that for any $\nu>0$ there exist positive constants C, C' such that
             \begin{equation}\label{eq2.6} 
             -H_p\lambda(x, \xi) \geqq CÊ\langle x\rangle^{-1- \nu}(\vert \xi\vert^2 +\vert x\vert^m)^ {\frac {1}{m}} - C' ,\quad  \forall (x, \xi) \in T^*(\R^n),
                \end{equation}
                (ii)\quad$\ds \lambda \in S(1,g),$ \\
              (iii) \ $\ds  [P, \lambda^w] - {\frac{1}{i}}(H_p \lambda)^w \in Op^wS_T(1,g).$
           \end{lem}
           Proof\\
        First of all on the support of $\chi(r)$  we have $\langle x\rangle^\frac {m}{2} \leqq \sqrt{p(x, \xi)} Ê\leqq C \vert \xi\vert.$ It follows that  $\vert \xi\vert \sim \langle \xi\rangle$ and\   $\vert \xi \vert  \leqq\vert \xi \vert +\langle x\rangle ^\frac{m}{2} \leqq C'\vertÊ\xi \vert.$ Now
          \begin{equation}\label{eq2.7} 
          -H_p \lambda = \sum_{j=1}^6 A_j
                 \end{equation}
                 where the  $A_j$'s   are defined below.\\
             1)   $A_1=\big(H_p \langle x\rangle^{-1} \big)p^{\frac{1}{m}-\frac{1}{2}}a\psi_0(\theta)\chi(r)$. Since on the support of $\psi_0(\theta)$ we have $\vert a\vertÊ\leqq 2 \varepsilon \langle x\rangle,$ it is easy to see that
              \begin{equation}\label{eq2.8}
              \vert A_1\vert \leqq C_1 \varepsilon \langle x\rangle ^{-1} \vert\xi\vert^\frac{2}{m}\big(1 - \psi(\vert \theta \vert)\big)\chi(r).
                \end{equation}
                2)  $A_2 =\langle x\rangle^{-1}p^{\frac{1}{m}-\frac{1}{2}}(H_pa)Ê\psi_0(\theta)\chi(r).$ By Lemma 2.1 $(i)$ we have
                  \begin{equation}\label{eq2.9}
                  A_2Ê\geqq C_2 \langle x\rangle^{-1}(\vert \xiÊ\vert + \langle x \rangle ^\frac {m}{2})^\frac {2}{m}
                  \big(1- \psi(\vert \thetaÊ\vert)\big) \chi(r) - C'_2.
                      \end{equation}
                      3)  $A_{3}=\langle
x\rangle^{-1}p^{\frac{1}{m}-\frac{1}{2}}a\psi_0'(\theta)(H_p\theta)
\chi(r)$. It follows from (\ref{eq2.3}), (\ref{eq2.4}) that
\begin{equation}\label{eq2.9b}
A_{3}=-p^{\frac{1}{m}-\frac{1}{2}}
|\theta|(H_p\theta)\psi'(|\theta|)\chi(r)
\end{equation}

4) $A_{4}=p^{\frac{1}{m}-\frac{1}{2}}(H_p\langle
a\rangle^{-\nu})\psi_1(\theta)\chi (r)$. Here we have $H_p\langle
a\rangle^{-\nu}=-\nu \langle a\rangle^{-2-\nu}aH_pa$. It follows from
(\ref{eq2.2}) that $A_{4}=\nu p^{\frac{1}{m}-\frac{1}{2}}|a|\langle
a\rangle^{-2-\nu} (H_pa)\psi(|\theta|)\chi(r)$.
Now on the support of $\psi(|\theta|)$ we have $\VE\langle
x\rangle\leq|a|$ and since $a\in S(\langle x\rangle,g)$ we have
$|a|\leq C\langle x\rangle$.
It follows from Lemma 2.1 $(i)$ that 
\begin{equation}\label{eq2.10}
A_{4}\geq C_3\langle x \rangle^{-1-\nu}(|\xi|+\langle x
\rangle^\frac{m}{2} )^\frac{2}{m}\psi(|\theta|)\chi(r) - C'_3.
\end{equation}

5) $A_{5}=-p^{\frac{1}{m}-\frac{1}{2}}(M_0-\langle
a\rangle^{-\nu})(H_p\theta)\psi_1'(\theta)\chi(r)$. It follows from
(\ref{eq2.3}) that 
\begin{equation}\label{eq2.11}
A_{5}=p^{\frac{1}{m}-\frac{1}{2}}(M_0-\langle
a\rangle^{-\nu})(H_p\theta)\psi'(|\theta|)\chi(r)
\end{equation}
We deduce from (\ref{eq2.9b}) and (\ref{eq2.11}) that
\begin{equation*}
A_{3}+A_{5}=p^{\frac{1}{m}-\frac{1}{2}}(M_0-\langle
a\rangle^{-\nu}-|\theta|)(H_p\theta)\psi'(|\theta|)\chi(r)
\end{equation*}
Now $H_p\theta=\langle x\rangle^{-1}H_pa+aH_p\langle
x\rangle^{-1}$. Since $|a|\leq 2\VE|\theta|$ 
on the support of
$\psi'(|\theta|)$ 
we deduce that $H_p\theta\geq C_4\langle x\rangle^{-1}|\xi|-C_5\geq -C_5$.
Taking $M_0\geq 2 $ and using the facts that $\psi'\geq 0$, $\chi\geq
0$ and $\VE\leq|\theta|\leq 2\VE$
on the support of
$\psi'(|\theta|)$
we obtain
\begin{equation}\label{eq2.12}
A_{3}+A_{5}\geq-C_6
\end{equation}

6) $A_{6} = (\langle x\rangle^{-1}a\psi_0(\theta)-(M_0-\langle
a\rangle^{-\nu})
\psi_1(\theta))p^{\frac{1}{m}-\frac{1}{2}}H_p[\chi(r)]$.
We have $\ds H_p[\chi(r)]=\frac{1}{\sqrt{p}}(H_p\langle
x\rangle^\frac{m}{2})\chi'(r)$. On the support of $\chi'(r)$ we have
$\langle x\rangle\sim|\xi|^\frac{2}{m}$;  this implies that 
\begin{equation*}
p^{\frac{1}{m}-\frac{1}{2}}|H_p[\chi(r)]|\leq
C|\xi|^{\frac{2}{m}-1}\frac{|\xi|\langle
  x\rangle^{\frac{m}{2}-1}}{|\xi|}|\chi'(r)|\leq C_7.
\end{equation*}
Therefore we obtain
\begin{equation}\label{eq2.13}
|A_6|\leq C_8.
\end{equation}
Gathering the estimates obtained in (\ref{eq2.8}) to  (\ref{eq2.13})
we obtain 
%%%%%%%%%%%%%%%%%%%%%%%%%%%%%%%%%%%%%%%%%%%%%%%%%%%%%%%%%%%%%%%%%%%%
%                                                                  %
%        Debut de la page 7                                        %
%                                                                  %
%%%%%%%%%%%%%%%%%%%%%%%%%%%%%%%%%%%%%%%%%%%%%%%%%%%%%%%%%%%%%%%%%%%%
\begin{equation}\label{eq2.14}
-H_p\lambda\geq C_9\langle x\rangle^{-1-\nu}(|\xi|+\langle
 x\rangle^\frac{m}{2})^\frac{2}{m} \chi(r)-C_{10}.
\end{equation}
Now on the support of $1-\chi(r)$ we have $|\xi|\leq C_{11}\langle
x\rangle^\frac{m}{2} $ so $ \langle x\rangle^{-1-\nu}(|\xi|+\langle
 x\rangle^\frac{m}{2})^\frac{2}{m}\leq C_{12}$. Therefore writing
   $1=1-\chi+ \chi$ and using (\ref{eq2.14}) we obtain (\ref{eq2.6}).
%%%%%%%%%%%% ATTENTION REF 2.6

$(ii)$ We use the symbolic calculus in the classes $S(M,g)$. We have
   $\langle x\rangle^{-1}\in S(\langle x\rangle^{-1},g)$, $a\in
   S(\langle x\rangle ,g)$, $p\in  S(\langle \xi\rangle^{2},g)$
so $p^{\frac{1}{m}-\frac{1}{2}}\in S(\langle
   \xi\rangle^{\frac{2}{m}-1},g)$ since $p\geq C>0$ on $\supp
   \chi(r)$.
Moreover $\chi(r)\in S(1,g)$ and on $\supp\chi(r)$ we have $\langle
   x\rangle^\frac{m}{2}\leq C|\xi|$.
It follows that $\lambda\in S(\langle
   \xi\rangle^{\frac{2}{m}-1},g)\subset S(1,g)$.

$(iii)$ By the symbolic calculus $\{\lambda,V\}\in
   S_T(\langle \xi\rangle^{\frac{2}{m}-1}\langle x\rangle^m\langle
   x\rangle^{-1}\langle \xi\rangle^{-1}  ,g) $. 
Since we have $\langle
   x\rangle^\frac{m}{2}\leq C|\xi|$ on its support we will have $\langle
   x\rangle^{m-1} \langle \xi\rangle^{\frac{2}{m}-2}\leq
   C|\xi|^{\frac{2}{m}(m-1)} \langle \xi\rangle^{\frac{2}{m}-2}\leq
   C'$.
Therefore $\{\lambda,V\}\in S_T(1,g)$. Now if $b\in  S_T(\langle
   x\rangle^\frac{m}{2}  ,g)$ we have 
$\{\lambda,b\xi_j\}\in S(\langle
   \xi\rangle^{\frac{2}{m}-1}\langle x\rangle^\frac{m}{2}|\xi|\langle 
   x\rangle^{-1}\langle \xi\rangle^{-1},g)$ and since $\langle
   x\rangle^\frac{m}{2}\leq C|\xi|$ we have $\langle
   x\rangle^{\frac{m}{2}-1}\langle
   \xi\rangle^{\frac{2}{m}-1}\leq
   C|\xi|^{\frac{2}{m}(\frac{m}{2}-1)}\langle
   \xi\rangle^{\frac{2}{m}-1}\leq C'$ so $\{\lambda, b\xi_j\}\in
   S_T(1,g)$. 

Finally $[ {\rm
   Op}^w(p),\lambda^w]-\frac{1}{i}(H_p\lambda)^w\in S(\langle
   \xi\rangle^2\langle \xi\rangle^{\frac{2}{m}-1}\langle x
   \rangle^{-2}\langle \xi\rangle ^{-2},g)\subset{\rm Op}^wS(1,g)$. \cqfd

%%%%%%%%%%%%%%%%%%
End of the proof of Theorem 1.1.

Since $\lambda \in S(1, g)$ we can set $\ds M=1+ \sup
_{(x,\xi)\in\R^{2n}}|\lambda (x,\xi)|$. Let us introduce
  $N(t)=((M+\lambda^w)u(t), u(t))_{L^2(\R^n)}$. Then there exist
absolute constants $C_1>0$, $C_2>0$ such that $C_1 \|
u(t)\|^2_{L^2}\leq N(u(t))\leq C_2  \| u(t)\|^2_{L^2}$. Now 
\begin{equation*}
\frac{d}{dt}N(t)=((M+\lambda^w)\frac{\partial u}{\partial t}(t),
u(t))_{L^2}
+((M+\lambda^w)u(t),\frac{\partial u}{\partial t}(t))_{L^2}
\end{equation*}
Since $\ds \frac{\partial u}{\partial t}(t)=-iPu(t)+if(t)$ and $P^\ast
=P$ we obtain
\begin{equation*}
\begin{aligned}
\frac{d}{dt}N(t)=&i([P,\lambda^w]u(t),u(t))_{L^2}-2\Im
((M+\lambda^w)f(t),u(t))_{L^2}\\
=&-((-H_p\lambda)^wu(t),u(t))_{L^2}-2\Im
((M+\lambda^w)f(t),u(t))_{L^2}+O(\| u(t) \|^2_{L^2})
\end{aligned}
\end{equation*}
By lemma 2.2 $(iii)$.

Now by Lemma 2.2 $(i)$ and the sharp G\aa rding inequality, we obtain 
\begin{equation}\label{eq2.15}
((-H_p\lambda)^wu(t),u(t))_{L^2}\geq C\| \langle
  x\rangle^{-\frac{1+\nu}{2}} E_\frac{1}{m}u(t)\|^2_{L^2}-C'\|
  u(t)\|^2_{L^2} 
\end{equation}
On the other hand we have for any $\VE >0$
%% Fin page 7
\begin{equation}\label{eq2.16}
|((M+\lambda^w)f(t), u(t))_{L^2}|\leq\VE\|\langle
 x\rangle^{-\frac{1+\nu}{2}}E_\frac{1}{m}u(t)\|^2_{L^2}+C_\VE\| \langle
 x\rangle ^\frac{1+\nu}{2} E_{-\frac{1}{m}}f(t)\|^2_{L^2}
\end{equation}
Using (\ref{eq2.15}) and (\ref{eq2.16}) with $\VE$ small enough, we
obtain
\begin{equation*}
\frac{d}{dt}N(t)\leq -C_1 \|\langle
 x\rangle^{-\frac{1+\nu}{2}}E_\frac{1}{m}u(t)\|^2_{L^2}+C_2\| \langle
 x\rangle ^\frac{1+\nu}{2} E_{-\frac{1}{m}}f(t)\|^2_{L^2}+C_3N(t)
\end{equation*}
Integrating this inequality between 0 and $t $ (in [0,T])  and using
Gronwall's inequality, we obtain the conclusion of Theorem 1.1. \cqfd

Proof of theorem 1.2.

Let $\chi\in C^\infty_0(\R^+)$, $\chi(t)=1$ if $t\in[0,1]$,
$\chi(t)=0$ if $t\geq 2$. Recall that according to (\ref{eq1.14}) we
have $p=|\xi |^2+q(x,\xi)$ where $\ds
q(x,\xi)=\sum^n_{j,k=1}b^{jk}(x)\xi_j\xi_k$ and $b^{jk}\in S(\langle
x\rangle^{-\sigma_0},g)$. Let us set 
\begin{equation}\label{eq2.17}
A_{jk}=\frac{x_j\xi_k-x_k\xi_j}{\langle \xi\rangle},\ 1\leq j,k\leq n
\end{equation}
Then we have the following result.
\begin{lem}
Let $a$ be defined in Lemma 2.1. One can find positive constants
$C_0$, $C_1$ and $C_2$ such that if we set 
\begin{equation}\label{eq2.18}
-\lambda=\frac{a}{\ds
  (1+a^{2}+\sum_{j,k=1}^nA_{jk}^2)^\frac{1}{2}}p^{\frac{1}{m}- 
\frac{1}{2}}\chi\left(\frac{\langle 
 x\rangle^\frac{m}{2}}{\sqrt{p(x,\xi)}}\right) 
\end{equation}
then
\begin{equation*}
\begin{aligned}
(i)&\ -H_p\lambda \geq C_0 \langle x\rangle ^{-3} (|\xi | +\langle
  x \rangle
  ^\frac{m}{2})^\frac{2}{m}\sum_{j,k=1}^nA_{jk}^2-C_1\langle
  x\rangle^{-1-\sigma_0}(|\xi|+\langle x\rangle
  ^\frac{m}{2})^\frac{2}{m}-C_2,\\ 
(ii)&\ \lambda\in S(\langle \xi\rangle ^{\frac{2}{m}-1},g),\\
(iii)&\ [P,\lambda^w]-\frac{1}{i}(H_p\lambda)^w\in {\rm Op}^wS_T(1,g).
\end{aligned}
\end{equation*}
\end{lem}
Proof

First of all we have
\begin{equation}\label{eq2.19}
|H_pA_{jk}(x,\xi)|\leq C_1\frac{|\xi|}{\langle x\rangle^{\sigma_0}},
\  1\leq j,k\leq n,\ (x,\xi)\in T^*(\R^n).
\end{equation}
Indeed we have $\{|\xi|^2,A_{jk}\}=0$ and $\ds |\{ q,A_{jk}\}|\leq
  C_2\frac{|\xi|}{\langle x\rangle^{\sigma_0}}$. 

Let us set
\begin{equation}\label{eq2.20}
D=1+a^2+\sum_{j,k=1}^nA_{jk}^2.
\end{equation}
We claim that on the support of $\chi(\langle
x\rangle^\frac{m}{2}p^{-\frac{1}{2}})$ we have
%%%%%%%%%%%%%%%%%%%%%%%%%%%%%%%%%%%%%%%%%%%%%%%%%%%%
%                                                  %
%   dŽbut de la page 9                            %
%                                                  %
%%%%%%%%%%%%%%%%%%%%%%%%%%%%%%%%%%%%%%%%%%%%%%%%%%%%
\begin{equation}\label{eq2.21}
C_3\langle x\rangle^2\leq D\leq C_4\langle x\rangle^2
\end{equation}
for some positive constants $C_3$ and $C_4$.

Indeed a straightforward computation shows that
\begin{equation*}
(x.\xi)^2+\sum_{j,k=1}^n(x_j\xi_k-x_k\xi_j)^2\geq |x|^2|\xi|^2.
\end{equation*}
Since by Lemma 2.1 we have $\ds a(x,\xi)=\frac{x.\xi}{\langle \xi\rangle}$
for $|x|\geq R_0\gg 1$ and $|\xi|\geq C_5>0$ on the support of $\chi$
we deduce that $D\geq C_6\langle x\rangle^2$ when $|x|\geq R_0$. When
$|x|\leq R_0$ we have $\ds D\geq 1\geq \frac{1}{1+R_0^2}\langle x\rangle^2$.

Now we can write with $r(x,\xi)=\langle
x\rangle^\frac{m}{2}p^{-\frac{1}{2}}$, 
\begin{equation}\label{eq2.22}\left\{
\begin{array}{l} 
-H_p\lambda=I_1+I_2\\
  I_1=D^{-\frac{3}{2}}(D(H_pa)-\frac{1}{2}a(H_pD))p^{\frac{1}{m}-
  \frac{1}{2}}\chi(r)\\ 
 I_2= p^{\frac{1}{m}-\frac{1}{2}}\ aD^ {-\frac {1}{2}}H_p(\chi(r)) 
 \end{array}\right.
\end{equation}
We have 
\begin{equation*}
\begin{aligned}
DH_pa-\frac{1}{2}
a(H_pD)=&(1+\sum_{j,k=1}^nA_{jk}^2)H_pa+a^2H_pa-\frac{1}{2}
a(2aH_pa+2\sum_{j,k=1}^nA_{jk}
H_pA_{jk})\\
=&(1+\sum_{j,k=1}^nA_{jk}^2)H_pa-a \sum_{j,k=1}^nA_{jk}
H_pA_{jk}.
\end{aligned}
\end{equation*}
Using (\ref{eq2.17}) and (\ref{eq2.19}) we see that,
\begin{equation}\label{eq2.23}
|a|\sum_{j,k=1}^n|A_{jk}||H_pA_{jk}|\leq C_7|x|^2\frac{|\xi|}{\langle
 x\rangle^{\sigma_0}}.
\end{equation}
Morever by Lemma 2.1 we have on the support of $\chi (r)$, 
\begin{equation}\label{eq2.24}
p^{\frac{1}{m}-\frac{1}{2}} (1+\sum_{j,k=1}^nA_{jk}^2)H_pa\geq
  (1+\sum_{j,k=1}^nA_{jk}^2)(C_8(|\xi|+\langle
  x\rangle^\frac{m}{2})^\frac{2}{m} -C_9).
\end{equation}
Therefore (\ref{eq2.20}), (\ref{eq2.22}), (\ref{eq2.23}),
(\ref{eq2.24}) show that,
\begin{equation*}
I_1\geqq\big[C_{10}\langle x\rangle^{-3}(|\xi|+\langle
  x\rangle^\frac{m}{2})^\frac{2}{m}
  \sum_{j,k=1}^nA_{jk}^2-C_{11}\frac{|\xi|}{\langle
  x\rangle^{1+\sigma_0}}\big]  \chi(r).
\end{equation*}
On the support of $1-\chi(r)$ we have $|\xi|\leq\langle
x\rangle^\frac{m}{2}$ so we obtain,
\begin{equation}\label{eq2.25}
I_1\geq C_{12}\langle x\rangle^{-3}(|\xi|+\langle
  x\rangle^\frac{m}{2})^\frac{2}{m}\sum_{j,k=1}^nA_{jk}^2-
  C_{13}\frac{(|\xi|+\langle  
  x\rangle^\frac{m}{2})^\frac{2}{m}}{\langle
  x\rangle^{1+\sigma_0}}-C_{14}.
\end{equation}
On the other hand we have,
\begin{equation*}
|H_p(\chi(r)|=|p^{-\frac{1}{2}}\chi'(r)H_p\langle
 x\rangle^\frac{m}{2}|\leq\frac{C_{15}}{|\xi|}|\chi'(r)||\xi|\langle
 x\rangle^{\frac{m}{2}-1} .
\end{equation*}
It follows from (\ref{eq2.21}) and the estimate $|a|\leq C_{16}
\langle x\rangle$ that,
%  debut page 10
\begin{equation}\label{eq2.26}
|I_2|\leq C_{17},
\end{equation}
since $\langle x\rangle^{\frac{m}{2}-1}|\xi|^{\frac{2}{m}-1}\leq C_{18}$.

Then $(i)$ in lemma 2.3 follows from (\ref{eq2.22}), (\ref{eq2.25}) and
(\ref{eq2.26}). The proofs of $(ii)$ and $(iii)$ are the same as those in the proof
of lemma 2.2. \cqfd

End of the proof of Theorem 1.2.

We introduce as before, for $t$ in $(0,T)$.
\begin{equation*}
N(t)=((M_0+\lambda^w)u(t),u(t))_{L^2}
\end{equation*}
Where $M_0$ is a large constant. Then $N(t)\sim \| u(t)\|^2_{L^2}$.

Now using the equation and Lemma 2.3 $(iii)$ we can write,
\begin{equation*}
\frac{d}{dt}N(t)=-((-H_p\lambda)^wu(t),u(t))_{L^2}-2\Im
((M_0+\lambda^w)f(t),u(t))_{L^2}+O(\| u(t) \|^2_{L^2})
\end{equation*}
Since by (\ref{eq1.13}) and (\ref{eq2.18}) we have $\langle
x\rangle^{-2}A^{2}_{jk}= \ell_{jk}^2$, Lemma 2.3 $(i)$ and the sharp G\aa
rding inequality ensure that
\begin{equation*}
\begin{split}
\frac{d}{dt}N(t)\leq 
& -C_1\sum_{j,k=1}^n\| \langle
  x\rangle^{-\frac{1}{2}}E_\frac{1}{m}
  \ell_{jk}^wu(t)\|_{L^2}^2+C_2\|\langle
  x\rangle^{-\frac{1+\sigma_0}{2}}E_\frac{1}{m}u(t)\|_{L^2}^2\\
&+\|\langle
  x\rangle
  ^\frac{1+\sigma_0}{2}E_{-\frac{1}{m}}f(t)\|_{L^2}^2
 +C_3N(t).
\end{split}
\end{equation*}
It follows that for $0<t<T$,
\begin{equation}\label{eq2.27}
\begin{split}
N(t)+C_1\int_0^t\sum_{j,k=1}^n\| \langle
  x\rangle^{-\frac{1}{2}}E_\frac{1}{m}
  \ell_{jk}^wu(s)\|_{L^2}^2ds\leq & N(0)+C_2\int_0^T\|\langle
  x\rangle^{-\frac{1+\sigma_0}{2}}E_\frac{1}{m}u(s)\|_{L^2}^2ds\\
&+\int_0^T\|\langle
  x\rangle
  ^\frac{1+\sigma_0}{2}E_{-\frac{1}{m}}f(s)\|_{L^2}^2ds
  +C_3\int_0^tN(s)ds. 
\end{split}
\end{equation}
Using Theorem 1.1 to bound the second term in the right hand side and
then using the Gronwall inequality we obtain
\begin{equation*}
N(t)\leq C(T)(\| u(0)\|_{L^2}^2+\int_0^T\|\langle
  x\rangle
  ^\frac{1+\sigma_0}{2}E_{-\frac{1}{m}}f(t)\|_{L^2}^2dt).
\end{equation*}
Using again the inequality (\ref{eq2.27}) we obtain the conclusion of
Theorem 1.2. The proof is complete.\cqfd

      L.Robbiano: Department of Mathematics,  University of Versailles, 45 Av. des Etats Unis 78035 Versailles France.  email: luc.robbiano@math.uvsq.fr\\
      C.Zuily: Department of Mathematics, University of Paris 11 91405 Orsay cedex France. email: claude.zuily@math.u-psud.fr
 
              \end{document}

    \subsection{}
 \begin{thm}
blablabla
\end{thm}

\begin{defi}
blablabla
\end{defi}

%meme punition pour les remarques, les lemmes, les corollaires etc.

%pour une Žquation qui va se numŽroter automatiquement

\begin{equation}
blablabla
\end{equation}

%pour une Žquation qui va se numŽroter et ˆ laquelle on pourra se rŽfŽrer grace au "label"

\begin{equation}\label{eq1.7}
blabla
 \end{equation}

%pour une Žquation sans numŽro

\begin{equation*}
blablabla
\end{equation*}

%autres exemples d'Žquations

\begin{equation}\label{eq }\left\{
\begin{array}{l}
equation avec une acolade et alignŽe ˆ gauche
\end{array}\right.
\end{equation}

\begin{equation*}
\begin{split}
blabla
 \end{split}
\end{equation*}

\begin{equation*}
\begin{aligned}
Tf(t)&=\int^b_a k(t,s)\,f(s)\,ds\\
Wf(t)&=\int^t_a K(t,s)\,f(s)\,ds\,.
\end{aligned}
\end{equation*}

\begin{itemize}
\item[(i)]
\end{itemize}

\end{document}